\documentclass[twoside,10pt,reqno]{amsart}
\usepackage{amsmath,amscd,amssymb,epsfig}
\def\wb{\overline}
\def\wt{\widetilde}

\def\End{\operatorname{End}}

\def\GL{\operatorname{GL}}
\def\Id{\operatorname{Id}}
\def\e{\operatorname{e}}
\def\L{{\Lambda}}
\def\M{{\mathcal{M}}}
\def\T{{\mathcal{T}}}
\def\E{{\mathcal{E}}}
\def\ZZ{\mathbb{Z}}  
\def\QQ{\mathbb{Q}}
\def\RR{\mathbb{R}}
\def\CC{\mathbb{C}}
\def\NN{\mathbb{N}}

\def\sll{\mathfrak{sl}}
\def\slto{\sll(2|1)}
\def\Usl{{{U_h}(\mathfrak{sl}(2|1)})}
\newcommand{\ch}{\operatorname{ch}}
\newcommand{\qn}[1]{{\left\{#1\right\}}}
\renewcommand{\S}[1]{{\mathcal{S}}_#1}
\newcommand{\qum}[1]{\widetilde{#1}}
\newcommand{\h}{\ensuremath{\mathfrak{h}}}
\newcommand{\g}{\ensuremath{\mathfrak{g}}}
\newcommand{\p}[1]{\ensuremath{\bar {#1}}}
\newcommand{\D}[1]{\ensuremath{\Delta({#1})} }
\newcommand{\HH}{\ensuremath{\boldsymbol H}}
\newcommand{\I}{\ensuremath{\boldsymbol I}}
\newcommand{\vz}{\ensuremath{\qum{V}_{0}}}
\newcommand{\twist}{\ensuremath{\theta}}
\newcommand{\twistvz}{\ensuremath{\theta_{\vz\otimes\vz}}}
\newcommand{\Conway}{\nabla}

\newcounter{bibcount}

\newtheorem{prop}{\bf Proposition}[section]
\newtheorem{defi}[prop]{\bf Definition}
\newtheorem{lem}[prop]{\bf Lemma}
\newtheorem{theo}{\bf Theorem}

\newtheorem{rk}[prop]{Remark}
\newcommand{\epsh}[2]
         {\begin{array}{c} \hspace{-1.3mm}
        \raisebox{-4pt}{\epsfig{figure=#1.eps,height=#2}}
        \hspace{-1.9mm}\end{array}}

\begin{document} 
\title[Multivariable link invariants and the Alexander
  polynomial]{Multivariable link invariants arising from $\slto$ and
  the Alexander polynomial} 
\author{Nathan Geer}
\address{School of Mathematics\\ 
Georgia Institute of Technology\\ 
Atlanta, GA 30332-0160}
\email{geer@math.gatech.edu}
\author{Bertrand Patureau-Mirand}
\address{L.M.A.M., Universit\'e de Bretagne-Sud, BP 573\\
F-56017 Vannes, France }
\email{bertrand.patureau@univ-ubs.fr}
\date{\today}

\begin{abstract}
In this paper we construct a multivariable link invariant arising from
the quantum group associated to the special linear Lie superalgebra
$\slto$.  The usual quantum group invariant of links associated to
(generic) representations of $\slto$ is trivial.  However, we modify
this construction and define a nontrivial link invariant.  This new
invariant can be thought of as a multivariable version of the
Links-Gould invariant.  We also show that after a variable reduction
our invariant specializes to the Conway potential function, which is a
refinement of the multivariable Alexander polynomial.
\end{abstract}

\maketitle
\setcounter{tocdepth}{1}

\section*{Introduction}

  There are deep connections between quantum group theory and
low-dimensional topology.  For example, every representation of a
semisimple Lie algebra gives rise to a quantum group invariant of
knots and, more generally, links.  It is well known that similar
invariants exist in the setting of Lie superalgebras.  Most nontrivial
invariants arising from Lie superalgebras are only invariants of long
knots or (1,1)-tangles.  This is true because, in many cases the
super-dimension of a finite-dimensional module over a Lie superalgebra
is zero.  This implies that the corresponding deformed module has
quantum dimension zero.  The standard quantum link invariant arising
from such a module is trivial.  For this reason it can be difficult to
construct non-trivial link invariants arising from Lie superalgebras.
In this paper we construct a multivariable link invariant arising from
representations over the quantum group associated to $\slto$.
  
  Our construction uses the Reshetikhin-Turaev quantum group invariant.
In particular, let $F$ be the usual functor from the category of
framed tangles colored by  topologically free $\Usl$-modules of finite 
rank to the category of
$\Usl$-modules (see \cite{Tu}).  In Section \ref{S:thF'}, we define a
map $d$ from the set of typical representations of $\Usl$ to the ring
$\CC[[h]]$.  If $T_{\lambda}$ is a framed $(1,1)$-tangle colored by
representations of $\Usl$ such that the open string is colored by the
deformed typical module $\qum{V}(\lambda)$ of weight $\lambda$, then
$F(T_{\lambda})=x Id_{\qum{V}(\lambda)}$ for some $x$ in $\CC[[h]]$.
We set $F'(T_{\lambda})=x.d(\lambda)$.
\begin{theo} \label{thF'}
  The map $F'$ induces a well defined invariant of framed links
  colored by at least one typical representation of $\Usl$.  In other
  words, if $L$ is a framed link colored by $\Usl$-modules at least
  one of which is typical and the closure of $T_{\lambda}$ is equal to
  $L$ then the map given by $L\mapsto F'(T_{\lambda})$ is a well
  defined framed link invariant.
\end{theo}
The proof of Theorem \ref{thF'} is given in Section \ref{S:thF'}.  We
denote the framed link invariant of Theorem \ref{thF'} by $F'$.  This
invariant can be thought of as a renormalization of the usual quantum
invariant.  Similar renormalization were considered by J. Murakami \cite{Mur93},
Kashaev \cite{Kv} and Degushi \cite{D}.  The construction here differs
from theirs as we work with ribbon categories whereas their proofs use
a Markov trace for the colored braid group.


Every irreducible topologically free representation of finite rank of
$\Usl$ has a highest weight $\lambda\in\Lambda\simeq \NN\times\CC$.
Thus the isomorphism classes of such representations are indexed by
the set $\NN\times\CC$.  In Section \ref{S:MVA} we will use $F'$ and
modules of highest weight of the form $(0,\alpha)$ to show that there
exists a generalized multivariable Alexander link invariant $M$.  
 In particular, we will prove the following
theorem.
\begin{theo}\label{T:Mpoly2}
  Let $L'$ be a framed oriented link with $n$ ordered components. Let
  $L$ be the non-framed link which underlies $L'$ ($L$ is still
  oriented and has the same ordering on its components).  There exists
  a generalized multivariable Alexander link invariant $M$ with the
  following properties.
  \begin{enumerate}
  \item \label{I:Mpoly1} If $n=1$ then $M(L)$ takes values in
  ${(q_{1}-q_{1}^{-1})^{-1}(q_{1}q-q_{1}^{-1}q^{-1})^{-1}}
  \ZZ[q^{\pm1},q_1^{\pm1}]$.
  \item If $n\geq 2$ then $M(L)$ takes values in
  $\ZZ[q^{\pm1},q_1^{\pm1},\ldots,q_n^{\pm1}]$.
  \item If $(\alpha_1,\ldots\alpha_n)\in(\CC\setminus\{0,-1\})^n$ and
  the $i$-th component of $L'$ is colored by the weight module of
  weight $(0,\alpha_{i})$ then
$$F'(L')={\e^{-\sum lk_{ij}(2\alpha_i\alpha_j+\alpha_i+\alpha_j)h/2}}
  {M(L)|_{q=\e^{h/2},q_i=\e^{\alpha_ih/2}}}$$
\end{enumerate} 
where $lk_{ij}$ is the linking number of the $i$-th and $j$-th
components of $L$.
 \end{theo} 
As a consequence of the last point in Theorem \ref{T:Mpoly2}, we get
that up to a change of variable $(p-1/p)(pq-1/(pq))M(q,p,p\ldots,p)$
is nothing but the Links-Gould invariant (\cite{LG}).  The Links-Gould
invariant is a two variable quantum group invariant arising from an
explicit one-parameter family of representations of the general linear
Lie superalgebra $\mathfrak{gl}(2|1)$.  The invariant $M$ can be
thought of as a multivariable version of the Links-Gould invariant
(see remark \ref{R:LG}).

In Section \ref{S:Skein}, using techniques similar to Viro in
\cite{Viro} we extend $F'$ to an invariant of colored oriented
framed trivalent graphs.  Using this extension, in Proposition
\ref{skein} we give a complete set of skein relations to compute this
invariant.  With the use of these skein relations we are able to show
that $M$ specialize to the Conway potential function.  (For a good
history and a nice geometric construction the Conway potential
function see \cite{C}.)  In particular, in Section \ref{S:Skein} we
prove the following theorem.
\begin{theo}\label{T:Conway}
Let $\Conway$ be the Conway potential function of a link.  Then
$$\Conway(L)(t_{1},..., t_{n})|_{t_{k}=q_{k}^{2}}\,=\,i
M(L)(q,q_{1},...,q_{n})|_{q=i}.$$
\end{theo}
This generalizes the results of \cite{WIL} which state 
that the two-variable Links-Gould invariant
dominates the one-variable Alexander polynomial.

In a subsequent paper the authors plan to use the framed link
invariant $F'$ to construct a quantum invariant of $3$-manifolds.
Because of the representation theory of $\slto$, this invariant would
be very different from $3$-manifold invariant arising from Lie
algebras.  We also plan to generalize the constructions of the
invariants $F'$ and $M$ to other Lie superalgebras.  The corresponding
$\sll$-invariants should specialize to the multivariable invariants of
\cite{ADO} and to Kashaev's quantum dilogarithm invariants of links
\cite{Kv}.
\subsection*{Acknowledgments} 
It is a pleasure to thank Christian Blanchet and Thang Le for
instructive conversations and useful suggestions. 
The authors also thank the referee for their very careful reading of
this manuscript and for their constructive remarks.
\section{Preliminaries}

In the section we review background material that will be used in the
following sections.

A \emph{super-space} is a $\ZZ_{2}$-graded vector space $V=V_{\p
0}\oplus V_{\p 1}$ over $\CC$.  We denote the parity of a homogeneous
element $x\in V$ by $\p x\in \ZZ_{2}$.  We say $x$ is even (odd) if
$x\in V_{\p 0}$ (resp. $x\in V_{\p 1}$).  A \emph{Lie superalgebra} is
a super-space $\g=\g_{\p 0} \oplus \g_{\p 1}$ with a super-bracket $[\:
, ] :\g^{\otimes 2} \rightarrow \g$ that preserves the
$\ZZ_{2}$-grading, is super-antisymmetric ($[x,y]=-(-1)^{\p x \p
y}[y,x]$), and satisfies the super-Jacobi identity (see \cite{K}).
Throughout, all modules will be $\ZZ_{2}$-graded modules (module
structures which preserve the $\ZZ_{2}$-grading, see \cite{K}).

\subsection{The Lie superalgebra $\slto$ and its weight modules}
In this subsection we define $\slto$ and discuss some properties of
$\slto$-modules.  Modules of $\slto$ are different in nature than
modules over semi-simple Lie algebras.  For example, the category of
finite-dimensional $\slto$-modules is not semi-simple.

Let $A=(a_{ij})$ be the $2\times 2$ matrix given by $a_{11}=2$,
$a_{12}=a_{21}=-1$ and $a_{22}=0$.
\begin{defi}
Let $\slto$ be the Lie superalgebra generated by $h_{i}$, $e_{i},$ and
$f_{i}$, $i=1,2$, where $h_{1}$, $h_{2}$, $e_{1}$ and $f_{1}$ are
even, $e_{2}$ and $f_{2}$ are odd, and the generators satisfy the
relations
\begin{align*}
\label{R:LieSuperalg}
   [h_{i}, h_{j}]&=0,   & [h_{i}, e_{j}] &=a_{ij}e_{j},   &
   [ h_{i},f_{j}] &=-a_{ij} f_{j} &   [e_{i},f_{j}] &=\delta_{ij}h_{i},
\end{align*}
\begin{align*}
 [e_{2},e_{2}] &=[f_{2},f_{2}] =0, & [e_{1},[e_{1},e_{2}]]
 &=[f_{1},[f_{1},f_{2}]]=0.
\end{align*}
\end{defi}
Set $\L=\NN\times\CC$.  For every pair $(a_{1},a_{2})\in\CC^{2}$ Kac
\cite{K} defined a highest weight $\slto$-module $V(a_{1},a_{2}) $
with a highest weight vector $v_{0}$ having the property that
$h_{i}.v_{0}=a_{i}v_{0}$ and $e_{i}v_{0}=0$.  We say $a=(a_{1},a_{2})$
is a weight and to simplify notation we will write $V(a)$ for
$V(a_{1},a_{2}) $.  Kac showed that finite-dimensional irreducible
$\slto$-modules are characterized up to isomorphism by the elements of
$\L$.  Moreover, the weight modules $V(a)$, $a\in\L$ are divided into
two classes: typical and atypical.

There are many equivalent definitions for a weight module to be
typical (see \cite{K78}).  In the interest of space, we will give the
following characterization of typical modules (which easily follows
from Theorem 1 in \cite{K78}).  A $\slto$ weight module
$V(a_{1},a_{2})$ is \emph{typical} if and only if $a_{1}+a_{2}+1\neq
0$ and $a_{2}\neq 0$.  If $V(a_{1},a_{2})$ is not typical we say it is
\emph{atypical}.  The following lemma, due to Kac \cite{K78}, is
useful when working with typical modules.

\begin{lem}\label{L:splits}
Typical module are projective and injective in the category of
finite-dimensional $\slto$-modules.  In particular, a typical $V(a)$
splits in any finite-dimensional representation.
\end{lem}

Let $\S1\subset \Lambda$ be the set of atypical weights.  Let $\S2
\subset \Lambda^{2}$ be the set
$$\S2:= (\S1\times \Lambda) \cup (\Lambda \times \S1) \cup \{(a,b)\in
\Lambda^{2}: V(a)\otimes V(b) \text{ is not semi-simple}\}.$$
 
The character of a typical module is well known (see \cite{K78}),
using these formulas and Lemma \ref{L:splits} we have the following
results about the tensor product of two modules.

\begin{lem}\label{L:TensorPab}  Let $a=(a_{1},a_{2})$ and
  $b=(b_{1},b_{2})$ be two weights such that $(a,b)\in \Lambda^2
  \setminus \S2$ and $b_{1}=0$.  If $a_{1}\neq 0$ then
  $V(a_{1},a_{2})\otimes V(0,b_{2})$ is isomorphic to
\begin{equation}
\label{E:tensor}
 V(a_{1}, a_{2}+b_{2}) \oplus V(a_{1}+1, a_{2}+b_{2}) \oplus
 V(a_{1}-1, a_{2}+ b_{2}+1) \oplus V(a_{1}, a_{2}+b_{2}+1).
\end{equation}
If $a_{1}=0$ then $V(0,a_{2})\otimes V(0,b_{2})$ is isomorphic to the
direct sum \eqref{E:tensor} without the module $V(a_{1}-1,a_{2}+b_{2}+1)$.
\end{lem}
\begin{proof}
Let $\h=<h_1,h_2>$ be the Cartan
subalgebra of $\sll(2|1)$.
The characters are elements of the group-ring $\ZZ\h^*$ of the space of
weights $\h^*$.  A basis of $\h^*$ is given by
$(\epsilon_1,\epsilon_2)$ with $\epsilon_1(h_1)=\epsilon_2(h_2)=1$,
$\epsilon_1(h_2)=0$ and $\epsilon_2(h_1)=-1$.  Hence, if
$(a_1,a_2)\in\Lambda$, the weight $w$ defined by $w(h_i)=a_i$ is
represented in $\ZZ\h^*$ by the element
$\e^{a_1\epsilon_1+a_2(\epsilon_1+\epsilon_2)}$.  In \cite{K78}, Kac
gives a general formula for the character of a typical module.  In our
context, the character of $V(a_{1},a_{2})$ is given by
\begin{equation}\label{Eq:Charac}
\ch(V(a_{1}, a_{2}))=(1+\e^{\epsilon_1})(1+\e^{\epsilon_2})
\e^{a_1\epsilon_1+a_2(\epsilon_1+\epsilon_2)}
\frac{1-r^{a_1+1}}{\hspace{-4ex}1-r}\text{ with }r=\e^{\epsilon_2-\epsilon_1}.
\end{equation}
Hence $\ch(V(a_1, a_{2})).\ch(V(0, b_{2}))=(1+\e^{\epsilon_1})(1+\e^{\epsilon_2})
\e^{a_1\epsilon_1+(a_2+b_2)(\epsilon_1+\epsilon_2)}X$ with 
\begin{align*}
X&=(1+\e^{\epsilon_1})(1+\e^{\epsilon_2})\frac{1-r^{a_1+1}}{\hspace{-4ex}1-r}\\
&=(1+\e^{\epsilon_1+\epsilon_2}+e^{\epsilon_1}(1+r))
\frac{1-r^{a_1+1}}{\hspace{-4ex}1-r}\\ 
&=\frac{1}{1-r}\left((1-r^{a_1+1})+ \e^{\epsilon_1+
  \epsilon_2}(1-r^{a_1+1})+ \e^{\epsilon_1}(1-r^{a_1+2})+
\e^{-\epsilon_1}\e^{\epsilon_1+\epsilon_2}(1-r^{a_1})\right)
\end{align*}
Then using equation \eqref{Eq:Charac}, we see that the product
$\ch(V(a_1, a_{2})).\ch(V(0, b_{2}))$ is the sum of the characters of
the typical representations of \eqref{E:tensor}.
\end{proof}
\begin{rk}\label{R:DegTensor}
Lemma \ref{L:splits} implies that if $(a,b)\in\Lambda^2$, the module
$V(a)\otimes V(b)$ is always the direct sum of typical modules direct
sum the semi-direct product of atypical modules.
\end{rk}

\subsection{The quantization $\Usl$}\label{SS:quant}

Let $h$ be an indeterminate.  Set $q=e^{h/2}$.  
We adopt the following notations:
$$q^z=\e^{zh/2}\quad\textrm{and}\quad\qn z=q^z-q^{-z}.$$
\begin{defi}[\cite{FLV,MS}]\label{D:Usl} Let $\Usl$ be the $\CC[[h]]$-Hopf
  superalgebra generated by the elements $h_{i},E_{i}$ and $ F_{i}, $
  $i = 1,2$, subject to the relations:
\begin{align*}
 [h_{i},h_{j}] &=0, & [h_{i},E_{j}]=&a_{ij}E_{j}, &
 [h_{i},F_{j}]=&-a_{ij}F_{j},
\end{align*}
\begin{align*}
 [E_{i},F_{j}]=&\delta_{i,j}\frac{q^{h_{i}}-q^{-h_{i}}}{q-q^{-1}},  &
 E_{2}^{2}=&F_{2}^{2}=0, 
 \end{align*}
\begin{align*}
  E_{1}^{2}E_{2}-(q+q^{-1})E_{1}E_{2}E_{1}+E_{2}E_{1}^{2}&=0 &
  F_{1}^{2}F_{2}-(q+q^{-1})F_{1}F_{2}F_{1}+F_{2}F_{1}^{2}&=0
\end{align*}
where $[,]$ is the super-commutator given by $[x,y]=xy-(-1)^{\p x \p
 y}yx$.  All generators are even except for $E_{2}$ and $F_{2}$ which
 are odd.  The coproduct, counit and antipode given by
\begin{align*}
\label{}
   \D {E_{i}}= & E_{i}\otimes 1+ q^{-h_{i}} \otimes E_{i}, &
   \epsilon(E_{i})= & 0 & S(E_{i})=&-q^{h_{i}}E_{i}\\
   \D {F_{i}} = & F_{i}\otimes q^{h_{i}}+ 1 \otimes F_{i}, &
   \epsilon(F_{i})= &0 & S(F_{i})=&-F_{i} q^{-h_{i}}\\
   \D {h_{i}} = & h_{i} \otimes 1 + 1\otimes h_{i}, & \epsilon(h_{i})
   = & 0 & S(h_{i})= &-h_{i}
 \end{align*}
\end{defi}

Let $E'=E_{1}E_{2}-q^{-1}E_{2}E_{1}$ and $F'=F_{2}F_{1}-qF_{1}F_{2}$.
Khoroshkin, Tolstoy \cite{KTol} and Yamane \cite{Yam94} showed that
$\Usl$ has an explicit $R$-Matrix $R$.  In particular, they showed
that $R=\check{R}K$ where
\begin{equation}
\label{E:Rcheck}
\check{R}=exp_{q}(\qn1 E_{1}\otimes F_{1})exp_{q}(-\qn1 E'\otimes
F')exp_{q}(-\qn1 E_{2}\otimes F_{2}), 
\end{equation}
\begin{equation}
\label{E:K}
K=q^{-h_{1}\otimes h_{2}-h_{2}\otimes h_{1} -2h_{2}\otimes h_{2}}
\end{equation}
and $exp_{q}(x):= \sum_{n=0}^{\infty}x^{n}/(n)_{q}!$,
$(n)_{q}!:=(1)_{q}(2)_{q}...(n)_{q}$ and $(k)_{q}:=(1-q^{k})/(1-q)$.

We say a $\Usl$-module $W$ is topologically free of finite rank if it
is isomorphic as a $\CC[[h]]$-module to $V[[h]]$, where $V$ is a
finite-dimensional $\slto$-module.  Let $\M$ the category of
topologically free of finite rank $\Usl$-modules.  A standard argument
shows that $\M$ is a ribbon category (for details see \cite{G04B}).
Let $V,W $ be objects of $\M$.  We denote the braiding and twist
morphisms of $\M$ as
\begin{align*}
c_{V,W}:&V\otimes W \rightarrow W \otimes V, &
\twist_{V}:&V\rightarrow V\
\end{align*} 
respectively.  We also denote the duality morphisms of $\M$ as
\begin{align*}
\label{}
  b_{V} :&\CC[[h]]\rightarrow V\otimes V^{*}, & d_{V}': & V\otimes
V^{*}\rightarrow \CC[[h]]
\end{align*} 

Let $\T=Rib_\M$ be the ribbon category of ribbon graphs over $\M$ in
the sense of Turaev (see \cite{Tu}).  The set of morphisms
$\T((V_1,\ldots,V_n),(W_1,\ldots,W_m))$ is a space of formal linear
combinations of ribbon graphs over $\M$.  Let $F$ be the usual ribbon
functor from $\T$ to $\M$ (see \cite{Tu}).
  
In \cite{G04A}, Geer constructs a specific isomorphism of topological algebras
$$\alpha: \Usl \rightarrow U(\slto)[[h]]$$
which induces a functor $\slto-Mod_{f}\rightarrow \M, V\mapsto \qum{V}$, where $\slto-Mod_{f}$ is the category of finite dimensional $\slto$-module.  It is shown in \cite{G05} that $\qum{V}(a)$ is a $\Usl$ weight module and a deformation of $V(a)$.   In other words, the characters of $V(a)$ and $\qum{V}(a)$ are equal and
the super-spaces $\qum{V}(a)$ is equal to $V(a)[[h]]$.  Thus,
the representation theory of $\Usl$ is parallel to that of the Lie
superalgebra $\slto$.
  
It is well known that the super-dimension of any typical
representation of $\slto$ is zero.  The discussion above implies that
the quantum dimension of any deformed typical representation over
$\Usl$ is zero.  It follows that the functor $F$ is zero on all closed
ribbon graph over with at least one color which is a deformed typical
module.  For this reason it can be difficult to construct nontrivial
link invariants from $\slto$.

\section{Proof of Theorem \ref{thF'}}\label{S:thF'}

In this section we prove Theorem  \ref{thF'}.  To this end, we compute
the value of the long Hopf link colored by finite-dimensional
representations.  We use this to define a function $d$ from the set of
weights $\Lambda\setminus \S1$ to $\CC[[h]]$.  Finally, we use $d$ to
show that $F$ induces a well defined invariant of links.

\begin{defi}\label{D:scalar}
If $T\in\T(\qum{V}(a),\qum{V}(a))$ where $V(a)$ is a
finite-dimensional irreducible $\slto$ weight module then
$F(T)=x.\Id_{\qum{V}(a)}\in\End_{\Usl}(\qum{V}(a))$ for an
$x\in\CC[[h]]$.  We define the bracket of $T$ to be $$<T>=x$$ For
example, if $V, V'$ are modules of $\M$ such that $V'$ is irreducible,
we define
$$S'(V,V')=\left< \epsh{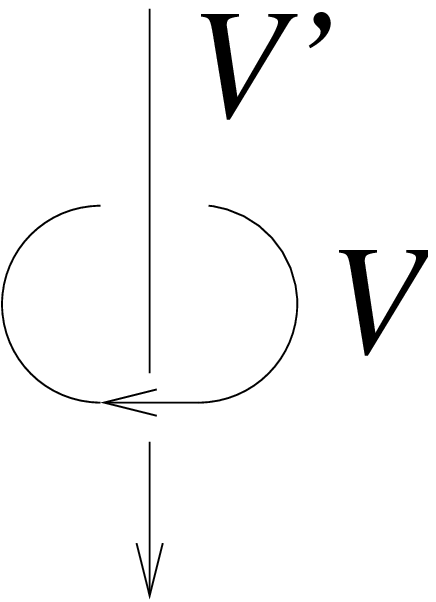}{8ex}\right>$$
\end{defi}
When $V=\qum{V}(a)$ and $V'=\qum{V}(b)$ are irreducible highest weight
modules with weights $a$ and $b$ in $\Lambda$ we write $S'(a,b)$ for
$S'(V,V')$.
\begin{prop}\label{P:S'} Let $a\in \Lambda\setminus \S1$ and $b\in
  \Lambda$ then
$$S'(a,b)= q^{-(2a_2+a_1+1)(2b_2+b_1+1)} 
  \frac{\qn{(a_1+1)(b_1+1)}}{\qn{b_1+1}} {\qn{b_2}}{\qn{b_2+b_1+1}}$$
\end{prop}
\begin{proof}
We will prove the proposition in three parts.  In the first part we
set up the notation, in the second we give two important facts and
then finally we give a calculation which completes the proof.

A direct calculation (using the character) shows that
$\qum{V}(a)=\qum{V}(a_{1},a_{2})$ has a basis of weight vectors
$v_{\p{j}, i}, w_{\p{j}, i}$ for $j=0,1$ and $i=0,...,a_{1}$ where
$v_{\p0,i}$, $w_{\p 0, i}$ are even and $v_{\p 1, i}, w_{\p1, i}$ are
odd.  The elements $h_{1}$ and $h_{2}$ act on these weight vectors as
follows:
\begin{align}
\label{E:actionhv}
   h_{1}v_{\p0, i} &= (a_{1}-2i) v_{\p0,i}, & h_{2}v_{\p0, i} &=
     (a_{2}+i) v_{\p 0,i}, \notag \\ 
   h_{1}v_{\p1, i} &= (a_{1}-2i+1) v_{\p1,i}, & h_{2}v_{\p 1, i} &=
     (a_{2}+i) v_{\p 1,i}, \notag\\ 
   h_{1}w_{\p0,i}&=(a_{1} -2i)w_{\p0,i}, & h_{2}w_{\p0,i}&=(a_{2}+
   i+1)w_{\p0,i}, \notag\\ 
   h_{1}w_{\p1,i}&=(a_{1} -2i-1)w_{\p1,i}, & h_{2}w_{\p1,i}&=(a_{2}+
   i+1)w_{\p1,i}.
\end{align}
Let $v_{b}$ be a highest weight vector of $\qum{V}(b)$.  Recall the
$R$-matrix is of the form $R=\check{R}K$.  From \eqref{E:actionhv} we
have
\begin{align}
\label{E:Kvv}
  K(v_{b}\otimes v_{\p0,i}) &= q^{-b_{1}(a_{2}+i) -b_{2}(a_{1}-2i)
  -2b_{2}(a_{2}+i)}v_{b}\otimes v_{\p0,i},\notag\\
  K(v_{b}\otimes v_{\p1,i}) &= q^{-b_{1}(a_{2}+i) -b_{2}(a_{1}-2i+1)
  -2b_{2}(a_{2}+i)}v_{b}\otimes v_{\p1,i},\notag\\
  K(v_{b}\otimes w_{\p0,i}) &= q^{-b_{1}(a_{2}+i+1) -b_{2}(a_{1}-2i)
  -2b_{2}(a_{2}+i+1)}v_{b}\otimes w_{\p0,i},\notag\\
  K(v_{b}\otimes w_{\p1,i}) &= q^{-b_{1}(a_{2}+i+1) -b_{2}(a_{1}-2i-1)
  -2b_{2}(a_{2}+i+1)}v_{b}\otimes w_{\p1,i}.
\end{align}
and $K$ is symmetric so acts by the same scalars on 
$v_{\p j,i}\otimes v_{b}$ and $w_{\p j,i}\otimes v_{b}$, respectively.\\

We now give two facts.  Let $v$ be any weight vector of $\qum{V}(a)$.\\
\begin{description}
  \item[Fact 1] $R(v_{b}\otimes v)=k(v_{b}\otimes v)$ where $k$ is the
  element of $\CC[[h]]$ given by the action of $K$ on $v_{b}\otimes
  v$.
  
  This fact follows from equation \eqref{E:Rcheck} and the property
  that $E_{i}v_{b}=0$ for $i=1,2$.
  \item[Fact 2] All the pure tensors of the element
  $(\check{R}-1)(v\otimes v_{b})\in \qum{V}(a)\otimes \qum{V}(b)$ are
  of the form $v' \otimes w'$ where $w'$ is a weight vector of
  $\qum{V}(b)$ and $v'$ is a weight vector of $\qum{V}(a)$ whose
  weight is of strictly higher order than that of $v$.
    
  Fact 2 is true because $E_{i}^{n}v$ (for $i=1,2$ and $n\in \NN$) is
  zero or a weight vector whose weight is of strictly higher order
  than the weight of $v$.
\end{description} 

We will now compute $S'(a,b)$ directly.  Let $V$ be an object of $\M$
and recall that the duality morphisms $b_{V} :\CC[[h]]\rightarrow
V\otimes V^{*}$ and $d_{V}': V\otimes V^{*}\rightarrow \CC[[h]]$ are
defined as follows.  The morphism $b_{V}$ is the $\CC[[h]]$-linear
extension of the coevaluation map on the underlying $\slto$-module.
In particular,
\begin{equation}
\label{E:bv}
b_{\qum{V}(a)}(1)=\sum_{j=0}^{1}\sum_{i=0}^{a_{1}}(v_{\p j,i}\otimes
v_{\p j,i}^{*}+w_{\p j,i}\otimes w_{\p j,i}^{*})
\end{equation}

As in the case of semi-simple Lie algebras we have
$$d_{\qum{V}(a)}' (v\otimes \alpha)= (-1)^{\p v \p
\alpha}\alpha(e^{h<\mu,\rho>}v)$$ where $v$ is a weight vector of
$\qum{V}(a)$ of weight $\mu\in\h^{*}$, $\rho$ is the half sum of the
positive even roots minus the half sum of the positive odd roots, and
$<,>$ is the natural non-degenerate bilinear form on $\h^{*}$ coming
from the super-trace.  From this a direct calculation shows
\begin{align}
\label{E:dva}
  d_{\qum{V}(a)}'(v_{\p j, i}\otimes v_{\p j,i}^{*}) &=
  (-1)^{\p j}e^{(-a_{2}-i)h}= (-1)^{\p j} q^{-2a_{2}-2i} \notag \\
  d_{\qum{V}(a)}'(w_{\p j, i}\otimes w_{\p j,i}^{*}) &=
  (-1)^{\p j}e^{(-a_{2}-i-1)h}=(-1)^{\p j} q^{-2a_{2}-2i-2}.
\end{align} 
Consider the element $S\in \End_{\Usl}(\qum{V}(b))$ given by
$$S=(Id_{\qum{V}(b)}\otimes
d_{\qum{V}(a)}')\circ(c_{\qum{V}(a),\qum{V}(b)}\otimes
Id_{\qum{V}(a)^{*}})\circ ( c_{\qum{V}(b),\qum{V}(a)}\otimes
Id_{\qum{V}(a)^{*}})\circ(Id_{\qum{V}(b)}\otimes b_{\qum{V}(a)}).$$ 
To simplify notation set $S=(X_{1})(X_{2})(X_{3})(X_{4})$ where
$X_{i}$ is the corresponding morphism in the above formula.  By
definition $S(v_{b})=S'(a,b)v_{b}$, so it suffices to compute
$S(v_{b})$.
\begin{align*}
S(v_{b})&=(X_{1})(X_{2})(X_{3})\Big(v_{b}\otimes
\sum_{j=0}^{1}\sum_{i=0}^{a_{1}}(v_{\p j,i}\otimes v_{\p
j,i}^{*}+w_{\p j,i}\otimes w_{\p j,i}^{*})\Big)\\
  &= (X_{1})(X_{2})\Big( \sum_{i = 1}^{a_{1}}\big(q^{-b_{1}(a_{2}+i)
  -b_{2}(a_{1}-2i) -2b_{2}(a_{2}+i)}v_{\p0,i}\otimes v_{b}\otimes
  v_{\p0,i}^{*} \\
  & \hspace{70pt} + q^{-b_{1}(a_{2}+i+1) -b_{2}(a_{1}-2i)
  -2b_{2}(a_{2}+i+1)}w_{\p0,i}\otimes v_{b}\otimes w_{\p0,i}^{*}\\
  &\hspace{70pt} +q^{-b_{1}(a_{2}+i) -b_{2}(a_{1}-2i+1)
  -2b_{2}(a_{2}+i)}v_{\p1,i}\otimes v_{b}\otimes v_{\p1,i}^{*}\\
  &\hspace{70pt} + q^{-b_{1}(a_{2}+i+1) -b_{2}(a_{1}-2i-1)
  -2b_{2}(a_{2}+i+1)}w_{\p1,i}\otimes v_{b}\otimes
  w_{\p1,i}^{*}\big)\Big)\\
  & =(X_{1})\Big( \sum_{i = 1}^{a_{1}}\big(q^{-2b_{1}(a_{2}+i)
  -2b_{2}(a_{1}-2i) -4b_{2}(a_{2}+i)} v_{b}\otimes v_{\p0,i}\otimes
  v_{\p0,i}^{*} \\
  & \hspace{70pt} + q^{-2b_{1}(a_{2}+i+1) -2b_{2}(a_{1}-2i)
  -4b_{2}(a_{2}+i+1)} v_{b}\otimes w_{\p0,i}\otimes w_{\p0,i}^{*}\\
  &\hspace{70pt} +q^{-2b_{1}(a_{2}+i) -2b_{2}(a_{1}-2i+1)
  -4b_{2}(a_{2}+i)} v_{b}\otimes v_{\p1,i}\otimes v_{\p1,i}^{*}\\
  &\hspace{70pt} + q^{-2b_{1}(a_{2}+i+1) -2b_{2}(a_{1}-2i-1)
  -4b_{2}(a_{2}+i+1)} v_{b}\otimes w_{\p1,i}\otimes
  w_{\p1,i}^{*}\\
 &\hspace{70pt} + \sum_{k} w'_{k}\otimes v'_{k}\otimes z_{k}\big)\Big)\\
  &=q^{-2b_{1}a_{2} -2b_{2}a_{1} -4b_{2}a_{2}-2a_{2}} \sum_{i =
  1}^{a_{1}}\Big(q^{-2b_{1}i +4b_{2}i -4b_{2}i-2i} + q^{-2b_{1}(i+1)
  +4b_{2}i -4b_{2}(i+1)-2i-2}\notag \\
  &\hspace{70pt} -q^{-2b_{1}i -2b_{2}(-2i+1) -4b_{2}i-2i} -
  q^{-2b_{1}(i+1) -2b_{2}(-2i-1) -4b_{2}(i+1)-2i-2} \Big)v_{b}\notag\\
  &= \big(q^{-(2a_2+a_1+1)(2b_2+b_1+1)}
  \frac{\qn{(a_1+1)(b_1+1)}}{\qn{b_1+1}}
  {\qn{b_2}}{\qn{b_2+b_1+1}}\big)v_{b}
\end{align*}
where $z_{k}= v_{\p j,i}^{*}$ or $w_{\p j,i}^{*}$ (for some $j=1,2$
and $i=0,...,a_{1}$), $w'_{k}$ is a weight vector of $\qum{V}(b)$ and
$v'_{k}$ is a weight vector of $\qum{V}(a)$.  From Fact 2 we have that
$z_{k}(v'_{k})=0$.  Moreover, the first equality of the above equation
follows from \eqref{E:bv}, the second from \eqref{E:Kvv} and Fact 1,
the third from \eqref{E:Kvv} and Fact 2, and finally the fourth from
\eqref{E:dva}.  The key observation in this proof is that Facts 1 and
2 imply that in the above computation the only contribution of the
action of the $R$-matrix comes from $K$.
\end{proof}
%
We are lead to the following definition.  
\begin{defi}
  If $a\in\Lambda\setminus \S1$, we set
$$d(a)=\frac{\qn{a_1+1}}{\qn1 \qn{a_2}\qn{a_2+a_1+1}}$$ so that one has
  the symmetry
 \begin{equation}
\label{E:sym}
  \forall (a,b)\in(\Lambda\setminus \S1)^2,\quad d(b)S'(a,b)=d(a)S'(b,a).
\end{equation}\end{defi}
\begin{lem}\label{L:twist}
Let $a=(a_{1},a_{2}) \in \Lambda\setminus \S1$.  Then the value of the
twist $\twist_{\qum{V}(a)}$ is $q^{-2a_{2}(1+a_{1}+a_{2})}$.  In other
words,
$$\left< \epsh{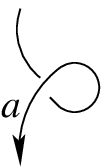}{6ex}\right> = q^{-2a_{2}(1+a_{1}+a_{2})}.$$
\end{lem}
\begin{proof}
The proof follow from \eqref{E:Kvv}, \eqref{E:dva} and Fact 1 in the
proof of Proposition \ref{P:S'}.
\end{proof}
\begin{lem}\label{L:Hopf2}
Let $c$ be the weight $(0,1)$.  Set $V_{0}:=V(c)$.  Then we have 
$$\left< \epsh{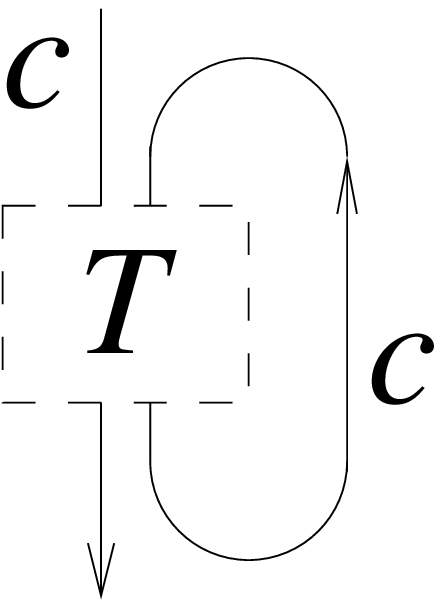}{8ex}\right>=\left< \epsh{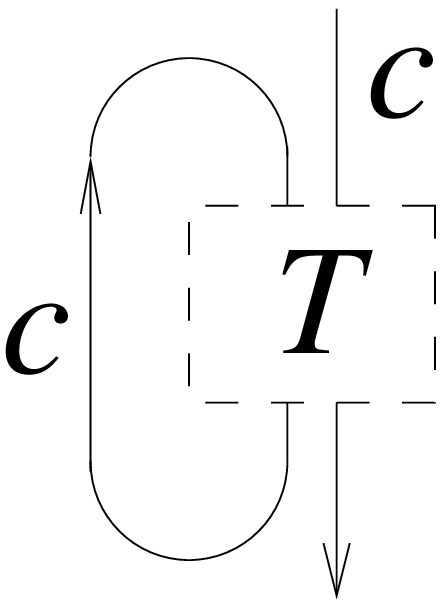}{8ex}\right>$$
for all $T\in\T\big((\vz, \vz), (\vz, \vz)\big)$. 
 \end{lem}
\begin{proof}
Set $E=\End_{\Usl}(\vz \otimes \vz)$.  Consider the following linear forms on
$E$:
$$tr_{L}(f)=(d_{\vz}\otimes Id_{\vz})\circ(Id_{\vz^{*}}\otimes
f)\circ(b'_{\vz}\otimes Id_{\vz}) \in \End_{\Usl}(\vz)\cong \CC[[h]],$$ 
$$tr_{R}(f)=(Id_{\vz}\otimes d'_{\vz}) \circ (f \otimes Id_{\vz^{*}})
\circ(Id_{\vz}\otimes b_{\vz}) \in \End_{\Usl}(\vz)\cong \CC[[h]].$$ 
To complete the proof it suffices to show that for any $T\in \T((\vz,
\vz),(\vz,\vz))$ we have \begin{equation}
\label{E:tr}
(tr_{L}\circ F)(T)=(tr_{R}\circ F)(T)
\end{equation} where $F$ is the functor described in subsection
\ref{SS:quant}.   

From Lemma \ref{L:TensorPab} we have that $\vz\otimes \vz\cong
\qum{V}(0,2)\oplus \qum{V}(1,2) \oplus \qum{V}(0,3)$.  Any element of
$E$ acts as a scalar on each summand and so the dimension of $E$ is 3.
We will now give a basis for $E$.

Consider the twist $\twistvz=F\left(\epsh{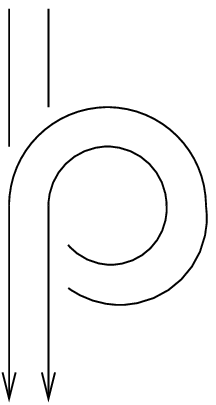}{6ex}\right)$ which is
an element of $E$.
From Lemma~\ref{L:twist} we have 
\begin{align*}
\label{}
  \twist_{\qum{V}(0,2)}=&q^{-12}  & \twist_{\qum{V}(1,2)}=&q^{-16}   &
  \twist_{\qum{V}(0,3)}=&q^{-24}.  
\end{align*}
Since 
$$\left|\begin{array}{ccc}
    1&q^{-12}&q^{-24}\\1&q^{-16}&q^{-32}\\1&q^{-24}&q^{-48}\end{array}\right|
    =(q^{-24}-q^{-16})(q^{-24}-q^{-12})(q^{-16}-q^{-12})\neq 0$$ 
we have that $(Id_{\vz \otimes \vz}, \twist_{\vz \otimes \vz},
\twist_{\vz \otimes \vz}^{2})$ form a basis of
 $E\otimes\CC[[h]][h^{-1}]$.  But $tr_{L}$ and $tr_{R}$ clearly have
the same values on this basis and thus equality \eqref{E:tr} holds as
$\CC[[h]]$ is an integral domain.
\end{proof}


\begin{lem}\label{key}
 Let  $(a,b)$ be any pair of weights belonging to
 $(\L\setminus\S1)^2$.  Then we have 
$$d(a)\left< \epsh{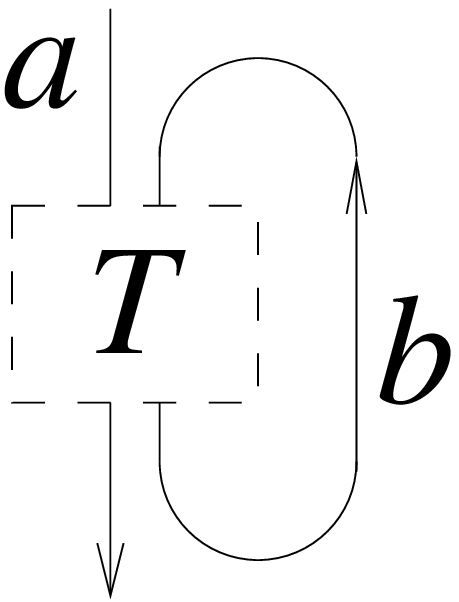}{8ex}\right>=d(b)\left< \epsh{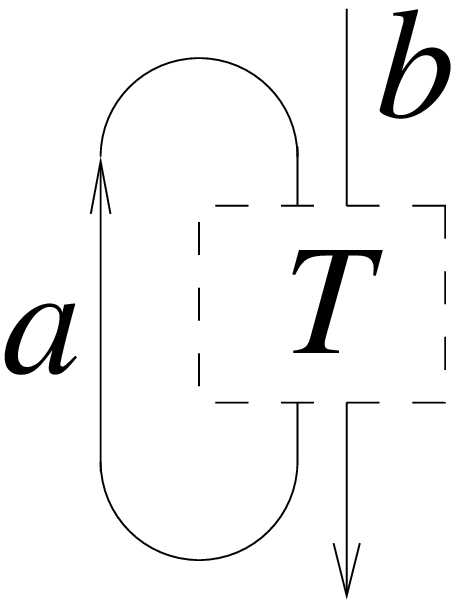}{8ex}\right>$$
for all $T\in\T\Big(\big(\qum{V}(a),\qum{V}(b)\big),
  \big(\qum{V}(a),\qum{V}(b)\big)\Big)$.
\end{lem}
\begin{proof}
Let $T\in\T\Big(\big(\qum{V}(a),\qum{V}(b)\big),
  \big(\qum{V}(a),\qum{V}(b)\big)\Big)$.  Let $c$ be the weight
  $(0,1)$.  By definition we have 
  \begin{align}
  \label{E:pict1}
 \left<\epsh{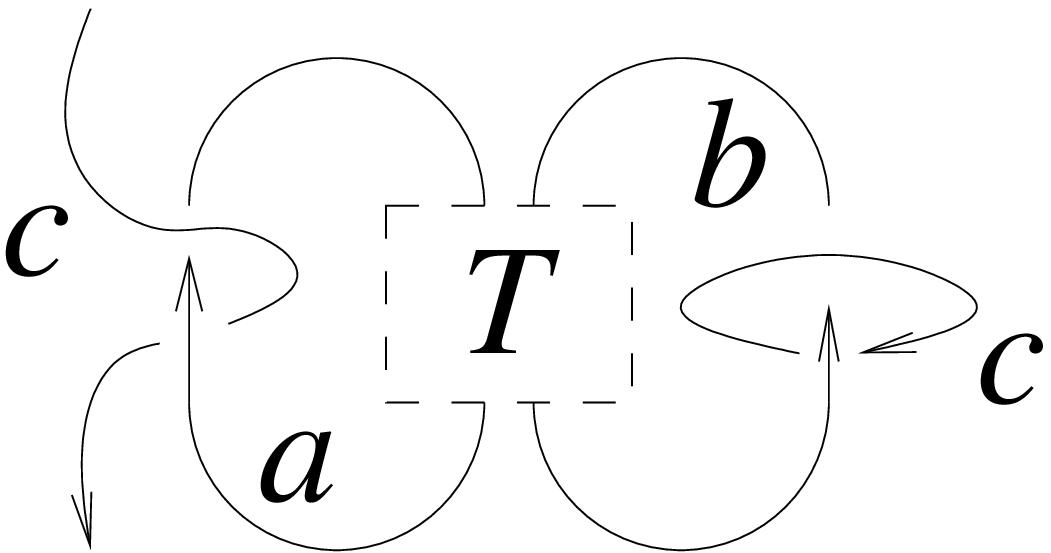}{8ex}\right>&=
  \left< \epsh{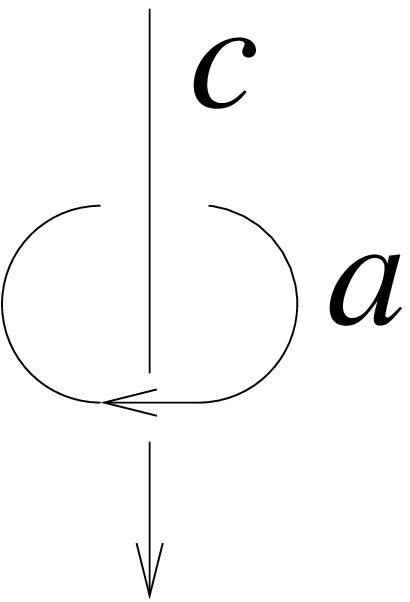}{8ex}\right>\left<\epsh{fig7}{8ex}\right>  \left<
  \epsh{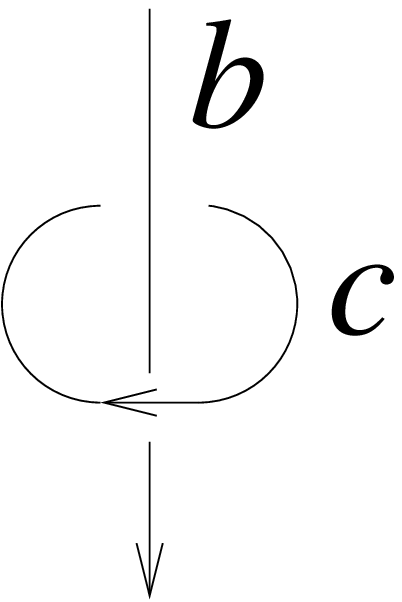}{8ex}\right> \notag\\  
  &=S'(a,c)S'(c,b)\left<\epsh{fig7}{8ex}\right>. 
  \end{align}

  Similarly,
  \begin{align}
  \label{E:pict2}
 \left<\epsh{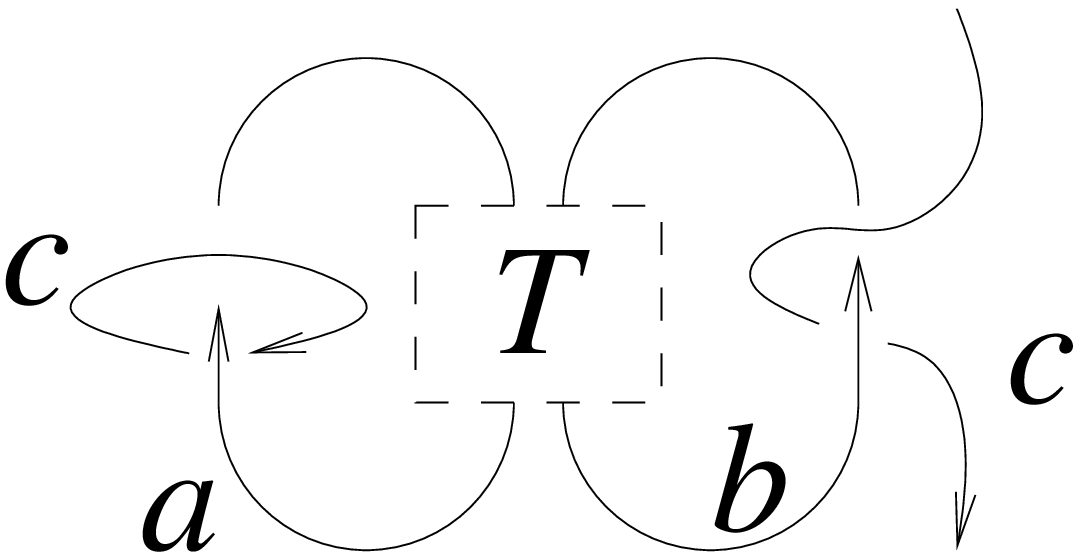}{8ex}\right>&=S'(b,c)S'(c,a)\left<\epsh{fig8}{8ex}\right>. 
  \end{align}

From Lemma \ref{L:Hopf2} we have that the left sides of equations
\eqref{E:pict1} and \eqref{E:pict2} are equal.  Thus, the results
follows from relation \eqref{E:sym}. 
\end{proof}
\begin{proof}[Proof of Theorem  \ref{thF'}]
Any closed ribbon graph $L\in\T(\emptyset,\emptyset)$ over $\M$ with
at least one edge colored by a typical module $\qum{V}(a)$ can be
represented as the closure of $T_a\in\T(\qum{V}(a),\qum{V}(a))$.  We
set $F'(L)=d(a)<T_a>$.  If $L$ can also be represented as the closure
of $T_b\in\T(\qum{V}(b),\qum{V}(b))$ for some typical weight $b$ then
there exits $T\in\T\Big(\big(\qum{V}(a),\qum{V}(b)\big),
\big(\qum{V}(a),\qum{V}(b)\big)\Big)$ such that
$T_a=\left<\epsh{fig7}{6ex}\right>$ and
$T_b=\left<\epsh{fig8}{6ex}\right>$ so by Lemma \ref{key} the
definition of $F'(L)$ does not depend on the choice of $T_a$.
\end{proof}

\section{The generalized multivariable Alexander invariant}\label{S:MVA}
Let $\rho_a:\Usl\rightarrow\End(V)$ be the representation
associated to the module $\qum V(a)\simeq V$.  When
$a_{2}\in\CC\setminus\{0,-1\}$ the representation $\rho_{(0,a_{2})}$
is four-dimensional.  We fix a super vector space $V$ of dimension
four with basis $B=(v_1,v_2,v_3,v_4)$ corresponding to the weight
vectors $(v_{\p0,0},w_{\p1,0},w_{\p0,0},v_{\p1,0})$ of the proof of
Proposition \ref{P:S'}. The super-grading of $V$ is given by
$\wb{v_i}=\wb{i+1}$.
\begin{lem}\label{L:mat} Up to equivalence, the matrices of the
  four-dimensional representation $\rho_{(0,a_{2})}$
  $(a_{2}\in\CC\setminus\{0,-1\})$ in the basis $B$ are:
{\small
$$
\begin{array}{ccc}
\left(\begin{array}{cccc} 
    0&0&0&0\\0&-1&0&0\\0&0&0&0\\0&0&0&1\end{array}\right)&
\left(\begin{array}{cccc} 
    0&0&0&0\\0&0&0&0\\0&0&0&0\\0&1&0&0\end{array}\right)&
\left(\begin{array}{cccc} 
    0&0&0&0\\0&0&0&1\\0&0&0&0\\0&0&0&0\end{array}\right)\\
\rho_{(0,a_{2})}(h_1)&\rho_{(0,a_{2})}(E_1)&\rho_{(0,a_{2})}(F_1)\\\\
\left(\begin{array}{cccc} 
    a_{2}&0&0&0\\0&a_{2}+1&0&0\\0&0&a_{2}+1&0\\0&0&0&a_{2}\end{array}\right)&
\left(\begin{array}{cccc} 
    0&0&0&q^{-a_{2}}\frac{\qn{a_{2}}}{\qn1}\\0&0&\frac{\qn{a_{2}+1}}{\qn1}&0\\
    0&0&0&0\\0&0&0&0\end{array}\right)&
\left(\begin{array}{cccc} 
    0&0&0&0\\0&0&0&0\\0&1&0&0\\q^{a_{2}}&0&0&0\end{array}\right)\\
\rho_{(0,a_{2})}(h_2)&\rho_{(0,a_{2})}(E_2)&\rho_{(0,a_{2})}(F_2)
\end{array}
$$
  }
\end{lem}
\begin{proof}
One can easily check that the relation of Definition \ref{D:Usl} are
satisfied and that the highest weight vector is $v_1$ whose weight is
$(0,a_{2})$.
\end{proof}
\begin{lem}\label{L:R-mat} There exists an
  $R(x,y,z)\in\GL(16;\ZZ[x^{\pm1},y^{\pm1},z^{\pm1}])$ such that, for
  all typical weights $a=(0,a_2)$ and $b=(0,b_2)$ the action of the
  $R$-matrix on $\qum{V}(a)\otimes \qum{V}(b)$ with respect to the
  basis $B\times B$ is given by $q^{-2a_2b_2}R(q,q^{a_2},q^{b_2})$.
\end{lem}
\begin{proof}
With the use of Lemma \ref{L:mat} and equations \eqref{E:Rcheck} and
\eqref{E:K} the lemma follows from a direct calculation. 
\end{proof}
\begin{proof}[Proof of Theorem \ref{T:Mpoly2}]
Choose $n$ complex numbers $\alpha_1,\ldots\alpha_n$ such that
$(1,\alpha_1,\ldots,\alpha_n)$ is a linearly independent family of the
$\QQ$-vector space $\CC$.  Then the ring map
$\phi:\ZZ[q^{\pm1},q_1^{\pm1},\ldots,q_n^{\pm1}]\rightarrow\CC[[h]]$
defined by $\phi(q)=\e^{h/2}$ and $\phi(q_i)=e^{\alpha_ih/2}$ is injective
since the family $\{\phi(q^{k_0}q_1^{k_1}\cdots q_n^{k_n}):
(k_0,\ldots ,k_n)\in\ZZ^{n+1}\}$ is free.  Let $L'$ be any framed link
with $n$ ordered components.  Color the $i$-th component of $L'$ with
the weight $(0,\alpha_{i})$.  Then from the definition of $F'$ and
Lemmas \ref{L:twist} and \ref{L:R-mat} we have
\begin{equation}
\label{E:F'}
F'(L')\in \frac{\e^{-\sum
lk_{ij}\alpha_i\alpha_jh}}{\qn{\alpha_k}\qn{\alpha_k+1}}
\operatorname{Im}(\phi)
\end{equation}
where $k\in\{1\cdots n\}$ is the number of
the component opened to compute $F'(L')$ (the denominator comes from
$d(0,\alpha_k)$).  

Let $L$ be the non-framed link which underlies $L'$.  If $n=1$ define
$$M(L):=(q_{1}-q_{1}^{-1})^{-1}(q_{1}q-q_{1}^{-1}q^{-1})^{-1}\phi^{-1}
\left(\e^{lk_{11}(\alpha_{1}^{2}+\alpha_{1})h}
\qn{\alpha_{1}}\qn{\alpha_{1}+1}F'(L')\right).$$ 
The first half of the theorem follow from this definition.  

Next we will show that if $n\geq 2$ then $F'(L')\in \e^{-\sum
lk_{ij}\alpha_i\alpha_jh}\operatorname{Im}(\phi)$.  For $i=1,2$ let
$T_{i}$ be a $(1,1)$-tangles whose closure is $L'$ and whose open
strand is the $i$-th component of $L'$.  From Theorem \ref{thF'} we
have $F'(T_{1})=F'(T_{2})$.  Then equation \eqref{E:F'} implies the
existence of Laurent polynomials $P_{1}$ and $P_{2}$ such that
 \begin{equation*}
{\e^{\sum lk_{ij}\alpha_i\alpha_jh}}F'(L') =d(0,\alpha_1)\phi(P_{1})
=d(0,\alpha_2)\phi(P_{2}).
\end{equation*}
It follows that
$(q_2-1/q_2)(q_2q-1/(q_2q))P_{1}=(q_1-1/q_1)(q_1q-1/(q_1q))P_{2}$.
Since $R=\ZZ[q^{\pm1},q_1^{\pm1},\ldots,q_n^{\pm1}]$ is an unique
factorization domain we have that $(q_1-1/q_1)(q_1q-1/(q_1q))$ divides
$P_{1}$.  Therefore we can conclude that
$${\e^{\sum lk_{ij}\alpha_i\alpha_jh}}F'(L')=
\phi\big(P_1/[(q_1-1/q_1)(q_1q-1/(q_1q))]\big).$$ 
For $n\geq 2$ we are now able to define 
$$M(L):=\phi^{-1} (\e^{\sum
  lk_{ij}(2\alpha_i\alpha_j+\alpha_i+\alpha_j)h/2}F'(L'))$$ 
where the additional correction is needed (see lemma \ref{L:twist}) to
make $M$ framing independent (i.e. a link invariant).

Because of Lemma \ref{L:R-mat}, $M(L)$ is independent of the choice
$\alpha=(\alpha_1,\ldots\alpha_n)$ lying in the dense subset of
$\CC^n$ defined by the condition: $(1,\alpha_1,\ldots\alpha_n)$ is a
linearly independent family of the $\QQ$-vector space $\CC$.  Now the
two maps $F'(L')$ and $\phi(M(L))$ depend continuously of $\alpha$ so
the relation between $F'$ and $M$ in Theorem  \ref{T:Mpoly2} is
valid for any $(\alpha_1,\ldots\alpha_n)\in(\CC\setminus\{0,-1\})^n$.
\end{proof}

\section{Skein relations}\label{S:Skein}

In this section we will give a complete set of skein relations for
$F'$ and use this to show that after a variable reduction $M$ is the
Conway potential function.  To this end, we extend $F$ and $F'$ to
invariants of colored oriented framed trivalent graphs.  For a
detailed account of similar extensions see the work of Viro
\cite{Viro}.

We will now define some normalized elementary $\Usl$-module morphisms.
Consider the following element of $V\otimes V$
$$\p b=v_1\otimes v_3 + q^{-1}v_3\otimes v_1 + v_4\otimes
v_2 - q^{-1}v_2\otimes v_4.$$
Remark that all typical representations $\qum V(0,a)$ have the same
underlying vector-space $V$ and the same basis $B$.
\begin{lem}\label{L:pbpd} The map $\p b_a : \CC[[h]]\rightarrow \qum V(0,a) \otimes \qum
  V_{(0,-1-a)}$ that sends $1$ to $q^{-a}\p b$ is $\Usl$-invariant
  and satisfies:
$$(\Id_{\qum V{(0,-1-a)}}\otimes\theta_{\qum V(0,a)})\circ c_{\qum V(0,a),
 \qum V{(0,-1-a)}}\circ\p b_a=\p b_{-1-a}$$
We use these morphisms to identify $\big(\qum V(0,a)\big)^*=\qum
V{(0,-1-a)}$. 
\end{lem}
With this identification, the evaluation map $d_{\qum
  V(0,a)}=\big(\qum V(0,a)\big)^*\otimes\qum V(0,a)\rightarrow \CC[[h]]$
  induces a family of bilinear maps: $\p d_a:\qum V(0,-1-a)\otimes\qum
  V(0,a)\rightarrow \CC[[h]]$ whose matrix in $B$ is given by
$$q^a\left( \begin {array}{rrrr}
 0&0&q&0\\0&0&0&1\\1&0&0&0\\0&-q&0&0\end {array} \right).$$ 
We also consider two families of morphisms (related to the so called
quantum Clebsch-Gordan coefficients):
$$\gamma_+^{a,b}:\qum V(0,a+b)\rightarrow\qum V(0,a)\otimes\qum V(0,b),\quad
\gamma_-^{a,b}:\qum V(0,a+b+1)\rightarrow\qum V(0,a)\otimes\qum
V(0,b)$$
given by:
\begin{align*}
\gamma_+^{a,b}(v_1)&=v_1\otimes v_1\\
\gamma_+^{a,b}(v_2)&=q^{-a}v_1\otimes v_2 + v_2\otimes v_1\\
\gamma_+^{a,b}(v_3)&=q^{-a}v_1\otimes v_3-q^b v_2\otimes v_4 + q^b
  v_3\otimes v_1 + q^{b+1} v_4 \otimes v_2 \\
\gamma_+^{a,b}(v_4)&=q^{-a}v_1\otimes v_4+v_4\otimes v_1\\
\\
\gamma_-^{a,b}(v_1)&=q^{b+1}\qn{a} \qn{a+1}v_1\otimes v_3 +q^b\qn{a+1}
 \qn{b+1} v_2 \otimes v_4\\ &+q^{-a-1}\qn{b} \qn{b+1} v_3\otimes v_1
 -q^{b+1}\qn{a+1} \qn{b+1} v_4\otimes v_2 \\
\gamma_-^{a,b}(v_2)&=\qn{a+1} \qn{a+b+1} v_2\otimes v_3
+q^{-a-1}\qn{b+1} \qn{a+b+1} v_3\otimes v_2\\
\gamma_-^{a,b}(v_3)&=\qn{a+b+1} \qn{a+b+2} v_3\otimes v_3\\
\gamma_-^{a,b}(v_4)&=q^{-a-1}\qn{b+1} \qn{a+b+1} v_3 \otimes v_4
+\qn{a+1} \qn{a+b+1} v_4\otimes v_3
\end{align*}
\begin{lem} \label{L:gamma} The map $\gamma_+^{a,b}$ (resp. $\gamma_-^{a,b}$) is
  $\Usl$-invariant and for $c=-1-a-b$ (resp. $c=-2-a-b$), it
  satisfies:
$$F\left(\epsh{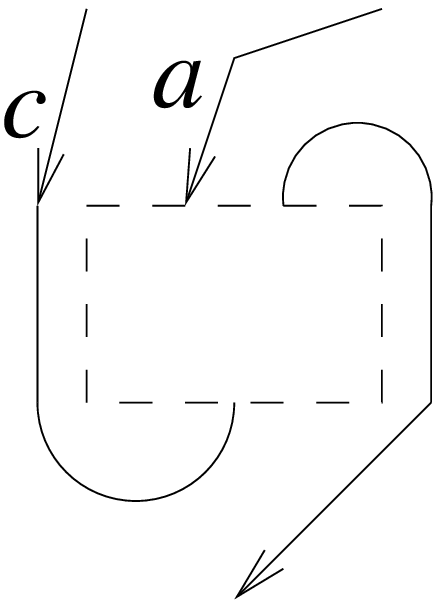}{10ex}\put(-22,0){\hbox{$\gamma_\epsilon^{a,b}$}}\right)=
  F\left(\epsh{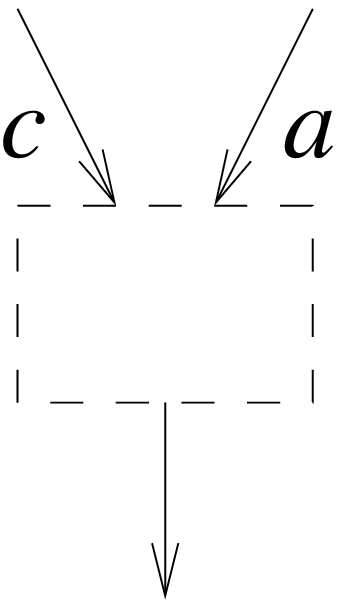}{10ex}\put(-20,0){\hbox{$\gamma_\epsilon^{c,a}$}}\right)$$ 
\end{lem}
This allow us to define $F'$ for a certain class of uni-trivalent ribbon graphs:
\begin{defi}
A trivalent tangle $T$ is a framed smooth embedding of a uni-trivalent
graph with oriented trivalent vertices in $\RR^3$.  Here the framing
is given by a continuous vector field along the image of $T$ nowhere
tangent to the image of $T$.  We impose that the three tangeant
vectors at a trivalent vertex are coplanar and two of them must not be
on the same half-line.  An orientation of a trivalent vertex $x$ is a
cyclic order on the three edges going to $x$. In the embedding $T$
this ordering must be given by the framing.
(At a trivalent vertex $x$ , the framing gives an orientation of the
affine plane tangeant to the image of $T$.  This orientation must be
compatible with the cyclic order of the tangeant of the three edges
going to $x$).

We denote by $\E(T)$ the set of oriented edges of $T$.
\begin{itemize}
\item $T$ is marked if it is given with a map from its trivalent
  vertices to the set $\{+,-\}\simeq\{\pm1\}$.
\item $T$ is colored if it is marked and if it is given with a map
  $f:\E(T)\rightarrow \CC$ such that for each oriented edge
  $\stackrel{\rightarrow}{e}\in\E(T)$, one has
  $f(\stackrel{\rightarrow}{e})=-1-f(\stackrel{\leftarrow}{e})$.
\item A coloring of $T$ is admissible if at each trivalent vertex
  marked with ``$+$'' (resp. ``$-$''), the sum of the values of $f$ at
  the three incoming edges is $-1$ (resp. $-2$).
\item $T$ is typical if it is colored with an admissible coloring for
  which $f$ takes values in $\CC\setminus\{0,-1\}$.
\item $T$ is closed if it has no univalent vertices.
\end{itemize}
\end{defi}
To any ribbon graph $G$ over $\M$ with edges colored by typical
modules $\qum V(0,a)$ and coupons marked by the morphisms
$\gamma_+^{a,b}$ and $\gamma_-^{a,b}$ one can associate a typical
trivalent tangle $T$. Lemma \ref{L:pbpd} and \ref{L:gamma} says that
$F(G)$ and $F'(G)$ depend only of $T$. We still denote by $F$ and $F'$
the induced maps on typical trivalent tangles.
\begin{rk}\label{R:FinZ} Let $T$ be a closed typical trivalent tangle.
We suppose $\dim(H_1(T))>1$ i.e. $T$ is not a knot.  Because of the
nice expression of $\p b$, $\p d$, $\gamma_+$, $\gamma_-$ and of the
$R$-matrix, we can make the following observation: As for links,
$F'(T)$ of a typical trivalent tangles has the form $q^{p_1}P_2$ where
$p_1$ is a degree $2$ polynomial in $\{f(e):e\in\E(T)\}$ and $P_2$ is
a Laurent polynomial in $\{q,q^{f(e)}:e\in\E(T)\}$. Unfortunately,
there is no canonical correction for the term $q^{p_1}$ and thus the
Laurent polynomial $P_2$ is only defined up to an unit.\\ For unclosed
typical trivalent tangles we have that the coordinates of $F(T)$ 
live in $\ZZ[q^{\pm1},q^{f(e)},q^{f(e)f(e')}:e,e'\in\E(T)]$.
\end{rk}
\begin{prop}\label{perturbation} The set of admissible colorings of a
  marked trivalent tangle $T$ is either empty or an affine space over
  $H_1(T,\CC)$. In the later case, the non typical colorings are a
  finite reunion of codimension $1$ subspaces. Furthermore, $F(T)$ is a
  continuous function of the coloring. Thus one can use this to compute
  $F(T)$ by perturbing the coloring: If $f$ is the coloring of $T$ and
  $c\in H_1(T,\CC)$, let $T_\epsilon$ be the marked trivalent tangle
  $T$ colored by $f+\epsilon c$ then
$$F\left(T\right)=\lim_{\epsilon\rightarrow0}F\left(T_\epsilon\right).$$
\end{prop}
\begin{proof}
Let $E$ be the set of oriented edges of $T$. The difference of two
admissible coloring of $T$ (as two elements of $\CC^E$) is then
equivalent to a $1$-cycle well defined up to a boundary (the relation
$f(\stackrel{\rightarrow}{e})=-1-f(\stackrel{\leftarrow}{e})$). 
The continuity is a consequence of Remark \ref{R:FinZ} (In fact, the
coordinates of $F\left(T_\epsilon\right)$ have the form
$q^{k\epsilon^2+\epsilon l(f)}P(q^\epsilon)$ where $k\in\ZZ$, $l$ is a
linear map of $\{f(e):e\in\E(T)\}$ and $P\in
\ZZ[q^{\pm1},q^{f(e)},q^{f(e)f(e')}:e,e'\in\E(T)][q^{\pm \epsilon}]$
is such that $P(1)$ is the corresponding coordinate of
$F\left(T\right)$).
\end{proof}
For $s,s'\in\{\pm1\}$, $a,b,c,d,e\in\CC$, satisfying
$a+b-e=-\frac{s+1}2$ and $c+d-e=\frac{s'-1}2$
(resp. $a-c+e=-\frac{s+1}2$ and $d-b+e=\frac{s'-1}2$) we set:
$$\I_{ab}^{cd}(ss')=\epsh{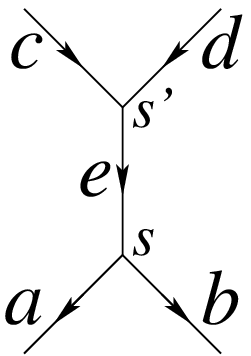}{10ex}\quad\text{and resp.}\quad 
\HH_{ab}^{cd}(ss')=\epsh{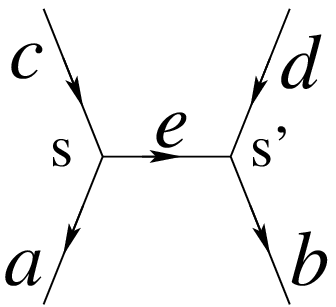}{9ex}$$
\begin{prop}\label{skein}$F$ and $F'$ satisfy the following skein
  relations. Furthermore up to the perturbation principle of
  proposition \ref{perturbation}, this set of skein relations is
  complete for closed typical trivalent tangles.  Indeed, we sketch in
  the proof an algorithm that allows to recursively compute the value
  of $F'$ on any closed typical trivalent tangle using these relations 
  and proposition \ref{perturbation}.
\begin{equation}\begin{array}{ll}\label{E:simple}
F\left(\epsh{fig24}{6ex}\right)=
q^{-2a(a+1)}F\left(\epsh{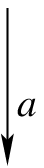}{6ex}\right)
&F\left(\epsh{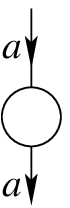}{6ex}\right)=
\qn{a}\qn{a+1}F\left(\epsh{fig25}{6ex}\right)\\
F\left(\epsh{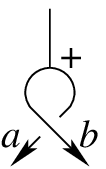}{6ex}\right)=
q^{2(a+1)(b+1)}F\left(\epsh{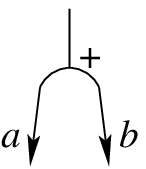}{6ex}\right)
&F\left(\epsh{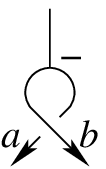}{6ex}\right)=
q^{2ab}F\left(\epsh{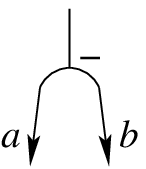}{6ex}\right)
\end{array}\end{equation}
\begin{align}
\label{E:IH1}F\left(\I_{ab}^{cd}(++)\right)&=
F\left(\HH_{ab}^{cd}(++)\right)\\ 
\label{E:IH2}\qn{c-a}\qn{c-a+1}F\left(\I_{ab}^{cd}(--)\right)&=
\qn{a+b}\qn{a+b+1}F\left(\HH_{ab}^{cd}(--)\right)\\
\label{E:IH3}i_1F\left(\I_{ab}^{cd}(-+)\right) +
i_2F\left(\I_{ab}^{cd}(+-)\right) &=h_1F\left(\HH_{ab}^{cd}(-+)\right)
+ h_2F\left(\HH_{ab}^{cd}(+-)\right)
\end{align}
\begin{equation}\label{E:Tab}\begin{array}{rl}
q^{2ab+a+b}F\left(\epsh{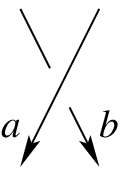}{6ex}\right)= &\frac1{\qn{a+b+1}} \left( 
\frac{q^{b}}{\qn{b}} F\left(\I_{ab}^{ba}(-+)\right) +
\frac{q^{-1-a}}{\qn{a+1}} F\left(\I_{ab}^{ba}(+-)\right) \right)\\
&-\frac1{\qn{b}\qn{a+1}}F\left(\HH_{ab}^{ba}(+-)\right)
\end{array}\end{equation}
\begin{equation}\label{E:Taa}\begin{array}{rl}
q^{2a(a+1)}F\left(\epsh{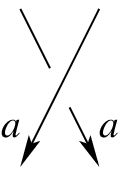}{6ex}\right)=
&-F\left(\epsh{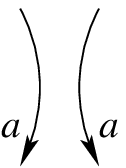}{6ex}\right)\\ 
&+\frac1{\qn{2a+1}} \left( \frac{q^{a}}{\qn{a}}
F\left(\I_{aa}^{aa}(-+)\right) + \frac{q^{-1-a}}{\qn{a+1}}
F\left(\I_{aa}^{aa}(+-)\right) \right)
\end{array}\end{equation}
\begin{equation}\label{E:split}
F'\left(\epsh{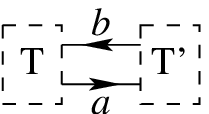}{6ex}\right)=\left\{\begin{array}{l}0\text{ if
}a\neq b\\ d(a)^{-1}F'\left(\epsh{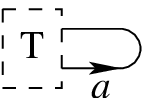}{5ex}\right)
F'\left(\epsh{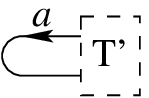}{5ex}\right)\text{ if }a=b\end{array}\right.
\end{equation}
where $
\begin{array}{l}
i_1= \qn{d-b}\qn{c+1}\qn{d+1} \\
i_2= -\qn{d-b}\qn{a}\qn{b}\\
h_1= -\qn{a+b+1}\qn{b}\qn{d+1}\\
h_2= \qn{a+b+1}\qn{a}\qn{c+1}.
\end{array}
$
\end{prop}
\begin{proof}
In general the existence of such relations is a basic consequence of the
representation theory of $\Usl$.  In this case the relations follow from lemma
\ref{L:TensorPab}.  We use a computer to find the coefficients.\\

We now give a sketch of an algorithm to compute $F'$ with these
relations:\\ 
Consider a regular planar projection of a closed typical trivalent
tangle.
\begin{enumerate}
\item One can convert it to a linear combination of planar graph using
relations (\ref{E:Tab}) and (\ref{E:Taa}).
\item Then the digons can be reduced by the second relation
(\ref{E:simple}).
\item We say that an edge $e$ is critical if its neighborhood is of
the form $\I_{ab}^{cd}$ with $c=a$ or $d=b$. We ignore for a moment
the problem of critical edges for which the relations (\ref{E:IH1}),
(\ref{E:IH2}) and (\ref{E:IH3}) can't be applied.
\item Consider the smallest $n$-gon. If it has two consecutive
vertices marked with the same sign, it can be reduced to a $(n-1)$-gon
using relation (\ref{E:IH1}) or (\ref{E:IH2}). then we reapply the
process from step 2.
\item If the smallest $n$-gon has its vertices alternatively marked
  $+$ and $-$ then use relation (\ref{E:IH3}) to obtain two
  $(n-1)$-gon and a $n$-gon as in step 4.
\item If during step 4 and 5 we deal with a critical edge then either
  it can be changed to a non critical edge using the perturbation
  principle, or it can be reduced by relation (\ref{E:split}).
\end{enumerate}
\end{proof}
In the rest of this section, we show that the specialization $q=i$ of
$M$ is essentially the Conway function:\\
For $\alpha\in\CC$, let $t^\alpha=q^{2\alpha}=\e^{\alpha h}$ and take
$h=i\pi$. Then for $k\in\ZZ$ we have 
$$q=i\qquad\qn{2k}=0\qquad\qn{2k+1}= (-1)^k2i\qquad\qn{\alpha+2}=
-\qn{\alpha}$$
$$\qquad\text{and}\qquad
\qn{\alpha}\qn{\alpha+1}=i\qn{2\alpha}=i(t^\alpha-t^{-\alpha}).$$

\begin{prop}\label{Askein} The specialization $h=i\pi$ of $F'$
  satisfies the following skein relations~:
\begin{equation}\label{E:Pab}
\qn{2(a+b)}F'\left(\epsh{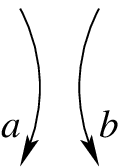}{6ex}\right) =
iF'\left(\I_{ab}^{ab}(+-)\right)-iF'\left(\I_{ab}^{ab}(-+)\right)
\end{equation}
\begin{equation}\begin{array}{ll}\label{E:Asimple}
F'\left(\epsh{fig24}{6ex}\right)=
t^{-a(a+1)}F'\left(\epsh{fig25}{6ex}\right)
&F'\left(\epsh{fig29}{6ex}\right)=
i\qn{2a}F'\left(\epsh{fig25}{6ex}\right)\\
F'\left(\epsh{fig30}{6ex}\right)=
t^{(a+1)(b+1)}F'\left(\epsh{fig31}{6ex}\right)
&F'\left(\epsh{fig32}{6ex}\right)=
t^{ab}F'\left(\epsh{fig33}{6ex}\right)
\end{array}\end{equation}
\begin{align}
\label{E:AIH1}F'\left(\I_{ab}^{cd}(++)\right)&=
F'\left(\HH_{ab}^{cd}(++)\right)\\ 
\label{E:AIH2}\qn{2(c-a)}F'\left(\I_{ab}^{cd}(--)\right)&=
\qn{2(a+b)}F'\left(\HH_{ab}^{cd}(--)\right)\\ 
\label{E:AIH3}\qn{2(a+b)}F'\left(\HH_{ab}^{cd}(+-)\right)
&= \qn{2d}F'\left(\I_{ab}^{cd}(-+)\right) -
\qn{2b}F'\left(\I_{ab}^{cd}(+-)\right)\\ 
\label{E:AIH4}\qn{2(a+b)}F'\left(\HH_{ab}^{cd}(-+)\right)
&=\qn{2a}F'\left(\I_{ab}^{cd}(+-)\right) -
\qn{2c}F'\left(\I_{ab}^{cd}(-+)\right)
\end{align}
\begin{equation}\label{E:ATab}\begin{array}{rl}
it^{ab}(t^{a+b}-t^{-a-b})F'\left(\epsh{fig26}{6ex}\right)=
F'\left(\I_{ab}^{ba}(-+)\right)
+t^{-a-b}F'\left(\I_{ab}^{ba}(+-)\right) \\
\end{array}\end{equation}
\begin{equation}\label{E:Asplit}
F'\left(\epsh{fig34}{6ex}\right)=\left\{\begin{array}{l}0\text{ if
}a\neq b\\ i\qn{2a}F'\left(\epsh{fig35}{5ex}\right)
F'\left(\epsh{fig36}{5ex}\right)\text{ if }a=b\end{array}\right.
\end{equation}
\end{prop}
\begin{proof}
Relations \eqref{E:Asimple}-\eqref{E:AIH2} and \eqref{E:Asplit} are
specializations of relations of Proposition \ref{skein} at $h=i\pi$.
Some new relations are only valid for $F'$ i.e. for closed trivalent
tangles: Relations \eqref{E:AIH3}, \eqref{E:AIH4} and \eqref{E:ATab}
are obtained from relation \eqref{E:Pab} by composing it with an $\HH$
or with the braiding.

We now prove the relation \eqref{E:Pab}:
Set $E=\T((\qum V(0,a),\qum V(0,b)),(\qum V(0,a),\qum V(0,b)))$ and
let $p$ be the element of $E$ given by
$$p=\qn{a+b}\qn{a+b+1}\qn{a+b+2}\Id -\qn{a+b}\I_{ab}^{ab}(+-)
-\qn{a+b+2}\I_{ab}^{ab}(-+)$$ 
Assume that $a+b\notin\ZZ$.  If $T$ is any element of $E$, denote
the closure (trace in $\T$) of $T$ by $\widehat T$.  Relation
\eqref{E:Pab} just says that for the specialization $h=i\pi$, one has
$F'(\widehat{T\circ p})_{|h=i\pi}=0$ for all $T\in E$.
Now $F(p)$ is just $\lambda p_2$ where
$\lambda=\qn{a+b}\qn{a+b+1}\qn{a+b+2}$ and $p_2$ is the projection on
the factor $\qum V(1,a+b)\subset\qum V(0,a)\otimes\qum V(0,b)$. This
is true because $\I_{ab}^{ab}(-+)$ (resp. $\I_{ab}^{ab}(+-)$) is
$\qn{a+b}\qn{a+b+1}$ times (resp. $\qn{a+b+1}\qn{a+b+2}$ times) the
projection $p_1$ (resp. $p_3$) on $\qum V(0,a+b)$ (resp $\qum
V(0,a+b+1)$) $\subset\qum V(0,a)\otimes\qum V(0,b)$.  We have
$F(T)=\sum \alpha_ip_i$ where the $\alpha_i$ are ``Laurent
polynomial'' (this is true because $p_j(v_i\otimes
v_i)=\delta_i^jv_i\otimes v_i$, thus $F(T)(v_i\otimes
v_i)=\alpha_iv_i\otimes v_i$ for $i=1,2,3$ and remark
\ref{R:FinZ}). Thus $F'(\widehat{T\circ p})=\alpha_2 F'(\widehat{p})$.
Now one can compute $F'(\widehat{p})= -\qn{a+b+2}-\qn{a+b}=
-(q+1/q)\qn{a+b+1}$ and $F'(\widehat{T\circ p})_{|h=i\pi}=0$.
\end{proof}
To prove Theorem \ref{T:Conway} we introduce a modified version of
Turaev's axioms for the Conway map (\cite{Tu86} section 4):
\begin{lem}\label{L:CAx}
The Conway function is the map uniquely determined by
\begin{enumerate}
\item $\Conway$ assigns to each ordered oriented link $L$ in $S^3$ an
  element of the field $\CC(t_1,\ldots,t_n)$ where $n$ is the number
  of components of $L$.
\item $\Conway(L)$ is unchange under (ambient) isotopy of the link
  $L$.
\item $\Conway(\text{unknot})=(t_1-t_1^{-1})^{-1}$.
\item \label{Ax:int} If $n\geq2$ then
$\Conway(L)\in\CC[t_1^{\pm1},\ldots,t_n^{\pm1}]$.
\item The one variable function
  $\wt\Conway(L)=\Conway(L)(t,t,\ldots,t)\in\CC[t^{\pm1}]$ is
  unchanged by a reordering of the components of $L$.
\item (Conway identity)
  $$\wt\Conway\left(\epsh{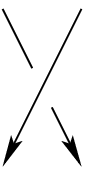}{4ex}\right)
  -\wt\Conway\left(\epsh{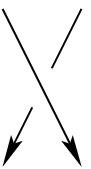}{4ex}\right) 
  =(t-t^{-1})\wt\Conway\left(\epsh{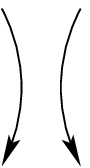}{4ex}\right)$$ 
\item \label{Ax:MDA}(Modified doubling axiom). If $L^+$ (resp. $L^-$) is obtained
  from the link $L=L_1\cup\cdots\cup L_n$ by replacing the component
  $L_i$ by its $(2,1)$-cable (resp. by its $(2,-1)$-cable) then
  \begin{align*}
    t_i\Conway(L^+)(t_1,\ldots,t_n)&-t_i^{-1}\Conway(L^-)(t_1,\ldots,t_n)=\\
    &\left(\prod_{j\neq i} t_j^{lk_{ij}}\right)(t_i^2-t_i^{-2})
    \Conway(L)(t_1,\cdots,t_{i-1},t_i^2,t_{i+1},\cdots,t_n)
  \end{align*}
\end{enumerate}
\end{lem}
\begin{proof}
The difference with the axioms given by Turaev (\cite{Tu86}) are :
\begin{itemize}
\item He only considers map with values in $\QQ(t_1,\ldots,t_n)$ and
  the axiom \ref{Ax:int} is replaced by
  $\Conway(L)\in\ZZ[t_1^{\pm1},\ldots,t_n^{\pm1}]$. 
  This change still
  allows us to consider $\wt\Conway$.
\item The axiom \ref{Ax:MDA} is replaced by the doubling axioms: If
  the link $L^+$ is obtained from the link $L=L_1\cup\cdots\cup L_n$
  by replacing the component $L_i$ by its $(2,1)$-cable then
  \begin{align*}
    \Conway(L^+)(t_1,\ldots,t_n)=(T+T^{-1})
    \Conway(L)(t_1,\cdots,t_{i-1},t_i^2,t_{i+1},\cdots,t_n) 
  \end{align*}
with $T=t_i\prod_{j\neq i} t_j^{lk_{ij}}$.
\end{itemize}
First the Conway function satisfies the modified doubling axiom
because
$$\Conway(L^-)(t_1,\ldots,t_n)=(T'+{T'}^{-1})
\Conway(L)(t_1,\cdots,t_{i-1},t_i^2,t_{i+1},\cdots,t_n)$$ with
$T'=t_i^{-1}\prod_{j\neq i} t_j^{lk_{ij}}$ (for example because $L^-$
is the miror image of the $(2,1)$-cable of the miror image of $L$ and
$\Conway(\text{miror}(L))=(-1)^{n+1}\Conway(L)$ (see \cite{Tu86})).

For the uniqueness, consider two maps $\Conway_1$ and $\Conway_2$ satisfying
axioms 1--7.  The proof of uniqueness given by Turaev has two steps: The
first step uses axioms 1--6 to show that the two one variable
specializations are the same: $\wt \Conway_1=\wt \Conway_2$.  This part
applies in our context without any change.

Now, using the modified doubling axiom, one can show by induction on
$\NN^n$ that for any $(a_1,\ldots,a_n)\in\NN^n$, one has
$\Conway_1(t^{2^{a_1}},\ldots,t^{2^{a_n}})=
\Conway_2(t^{2^{a_1}},\ldots,t^{2^{a_n}})$. But this implies that for
any fixed ordered link $L$ (with $n\geq2$ components), the Laurent
polynomial $\Conway_1(L)-\Conway_2(L)$ is zero.
\end{proof}
\begin{proof}[Proof of Theorem  \ref{T:Conway}] 
If $k\in\NN$, set $t_k=q_k^2$ and
$M'(q_1,\ldots,q_n)=i.M(i,q_1,\ldots,q_n)$.  It is clear that $M'$
satisfies the axioms 1--5 of lemma \ref{L:CAx} (we neglect the fact
that the values of $M'$ are a priori in the extension $\CC[q_k^{\pm1}]$ of
$\CC[t_k^{\pm1}]$).

Applying the braiding to equation \eqref{E:Pab} then using equation 
\eqref{E:Asimple} one can check that 
$$t^{a+a^2}F'\left(\epsh{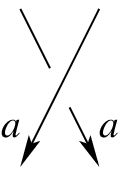}{4ex}\right)
-t^{-a-a^2}F'\left(\epsh{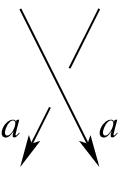}{4ex}\right)= 
(t^a-t^{-a})F'\left(\epsh{fig14}{4ex}\right).$$
This implies that $M'$ satisfies the axiom of the Conway
identity.

Next we show the modified doubling axiom holds.   From 
equation~\eqref{E:ATab} and its
mirror analog, one has 
\begin{equation}
\label{E:FI}
t^{2a+a^2}F'\left(\epsh{fig10}{4ex}\right)
-t^{-2a-a^2}F'\left(\epsh{fig11}{4ex}\right)=-i\I_{a,a}^{a,a}(-+).
\end{equation}
We choose a zero framing on $L$.  Hence the linking matrix $(lk_{ij})$ of $L$ satisfies the condition $lk_{ii}=0$.  
With this framing $L^{\pm}$ can
be obtain from $L$ by replacing $L_i$ by two parallel copies modified
in a small ball by a positive (or negative) crossing ($L^{\pm}$ inherits the framing of $L$ and its $i$\textsuperscript{th} component has framing $\pm1$). 
Combining the previous sentence with equations \eqref{E:Asimple} and \eqref{E:FI} we have 
$$t^{2a_i+a_i^2}F'(L^+)-t^{-2a_i-a_i^2}F'(L^-)=
i(t^{2a_i}-t^{-2a_i})F'(L)$$ where the $k$-th component of $L^+$,
$L^-$ and $L$ is colored by $a_k$ except the $i$-th component of $L$
which is colored by $2a_i$. Now as
$lk_{ik}(L^\pm)=lk_{ki}(L^\pm)=2lk_{ik}(L)$ for $k\neq i$ and
$lk_{ii}(L^\pm)=\pm1$, the framing correction gives 
\begin{align*}
t^{a_i}M'(L^+)&(t^{a_1},\ldots,t^{a_n})-t^{-a_i}M'(L^-)(t^{a_1},\ldots,t^{a_n})=\\
    &\left(\prod_{j\neq i} t^{lk_{ij}{a_j}}\right)(t^{2a_i}-t^{-2a_i})
    M'(L)(t^{a_1},\cdots,t^{a_{i-1}},t^{2a_i},t^{a_{i+1}},\cdots,t^{a_n})
\end{align*}
Hence $M'$ is the Conway function.
\end{proof}
Remark that Proposition \ref{Askein} is a complete set of skein
relations for a generalization of the Conway potential function
$\Conway$ to colored trivalent graphs.

\section{Some examples}
One can use the skein relations developed in Section \ref{S:Skein} to
compute the invariant $M$.  In this section we give the results of
such computations.   

For a knot $K$, the invariant $M$ is not a Laurent polynomial but it
has the form $M_0+P$ with
$M_0=\frac{1}{(q_1-q_1^{-1})(qq_1-(qq_1)^{-1})}=M(\text{unknot})$ and
$P\in\ZZ[q^\pm,q_1^\pm]$. In fact $M(K)$ is just the Links-Gould
invariant. More precisely we get 
$$LG_K(q,p)_{|p=q_1\sqrt{q}}=(q_1-q_1^{-1})(qq_1-(qq_1)^{-1})M(K)(q,q_1)$$
with the convention of \cite{Wit}. In this last paper, the values of
LG for the first prime knots (up to 10 crossings) are computed. The
corresponding values of $M$ can be deduced. As an example, in
\cite{Wit} the value for the trefoil ${3_1}$ is presented by
$$LG_{3_1}=1+2q^2-(q+q^3)(p^2+p^{-2})+q^2(p^4+p^{-4})$$
and our computation gives
$$M(3_1)=M_0\,+\,q^2(qq_1^2+(qq_1^2)^{-1})$$
this result and several other computations for knots are in agreement
with the computations of \cite{Wit}.

\begin{rk}\label{R:LG}
In \cite{WLK} the Links-Gould invariant is computed using an
$R$-matrix of $U_q{\mathfrak gl}(2|1)\simeq U_q{\mathbb T}_1\otimes
U_q\sll(2|1)$ (isomorphism of Hopf algebra) where $U_q{\mathbb T}_1$
is the (co-)commutative
 Hopf algebra of polynomials in one primitive
variable $c$.  
They consider the representation $\qum
V(0,\alpha)$ (for a generic value of $\alpha$) that is obtained from
the one of Lemma \ref{L:mat} by making $c$ acting by the scalar
$\alpha$ on it.
Hence, they compute the Reshetikhin-Turaev invariant
of a $(1,1)$-tangle $T$ where each component is colored by $\qum V(0,\alpha)$ with a
$R$-matrix that differs form ours by a scalar $q^{2\alpha(\alpha+1)}$.
This scalar exactly corrects the framing of the tangle.

The Links-Gould invariant is given by $<T>$ as in Definition
\ref{D:scalar}.  Remark that the specialization $b=a$ of Lemma
\ref{key} (which is not trivial for tangles with several components)
gives a proof that the Links-Gould invariant is a well defined
invariant of the link closure of $T$.

So again with the convention of \cite{Wit} where
$p=q^{\alpha+\frac12}$, we get that for any link $L$,
$$LG(L)(q,q^{\alpha+\frac12})=
(q^\alpha-q^{-\alpha})(q^{\alpha+1}-q^{-\alpha-1})
M(L)(q,q^\alpha,q^\alpha,\ldots,q^\alpha).$$ 
\end{rk}

It is more interesting to see the value for links with several
components.  For example, let $H$ be the Hopf link (with negative
crossings) then
$$M(H)(q,q_1,q_2)=q$$
and thus $M(H)$ does not depend of the two colors.  This result can be
deduced from Proposition \ref{P:S'} but we use it to illustrate the
skein relations:
$$\begin{array}{rll}
F'\left(\epsh{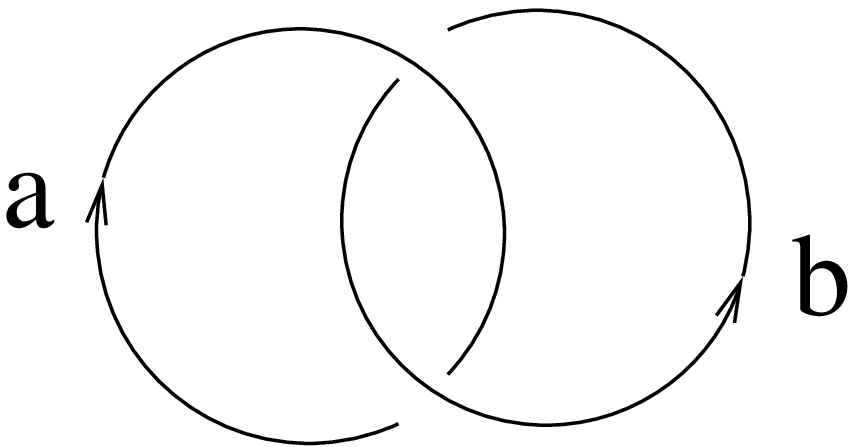}{5ex}\right)&
 =\frac{q^{2ab+a+b}q^{-a}}{\qn{a+b+1}\qn{a}}F'\left(\epsh{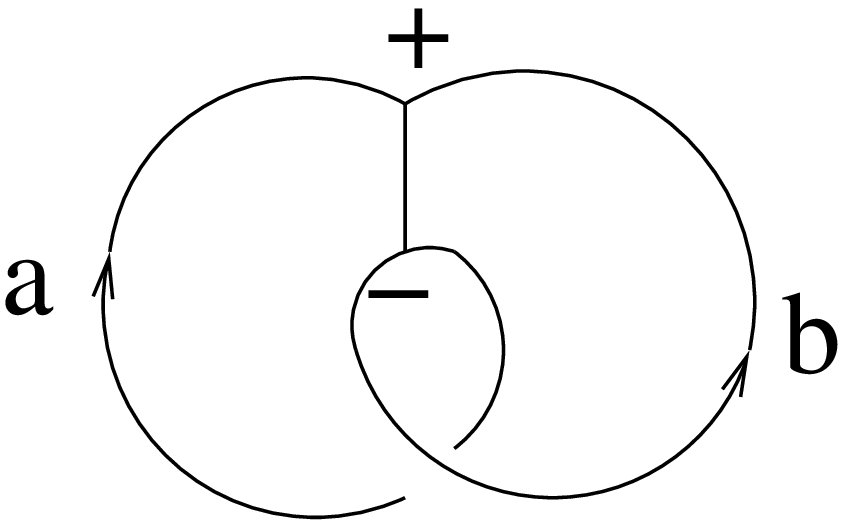}{5ex}\right)&
+\frac{q^{2ab+a+b}q^{b+1}}{\qn{a+b+1}\qn{b+1}}F'\left(\epsh{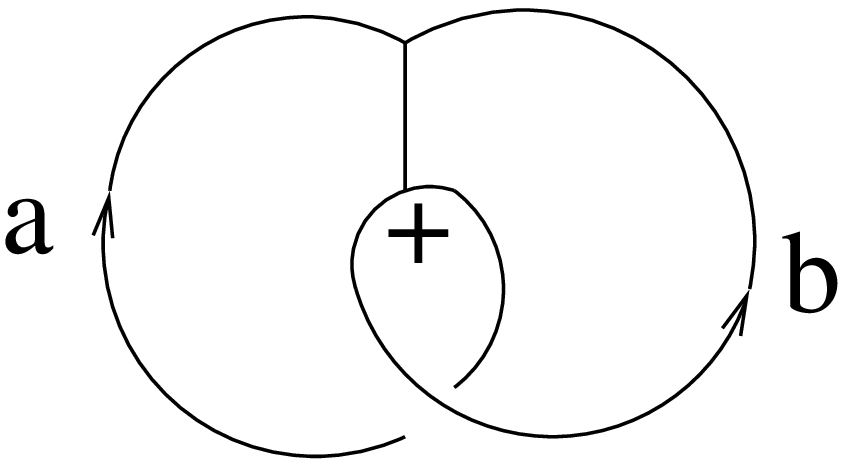}{5ex}\right)\\\\
&&-\frac{q^{2ab+a+b}}{\qn{a}\qn{b+1}}F'\left(\epsh{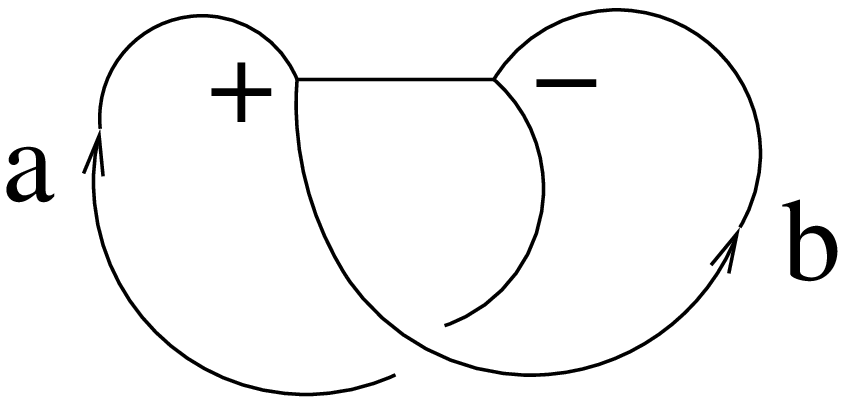}{5ex}\right) 
\end{array}$$
$$\begin{array}{l}=\frac{q^{2ab+a+b}q^{-a}}{\qn{a+b+1}\qn{a}}q^{2ab}+
\frac{q^{2ab+a+b}q^{b+1}}{\qn{a+b+1}\qn{b+1}}q^{2(a+1)(b+1)}-
\frac{q^{2ab+a+b}}{\qn{a}\qn{b+1}}q^{2a(b+1)}\\\\
=q^{4ab+2a+2b}q\quad\text{and }M(H)=q.
\end{array}$$
The first equality is obtained from the miror analog of \eqref{E:Tab}.
For the second we use equations \eqref{E:simple}.  (Remark that the
image by $F'$ of any typical planar ``Theta'' ($\Theta$) graph is
equal to $1$.)\\

The Borromean link $B$ ($B=L6a4$ in the
Thistlethwaite link table) has $3$ symmetric components:
\begin{eqnarray*}
M(B)(q,q_1,q_2,q_3)=\delta(q)+\delta(q_1)\delta(q_2)\delta(q_3)\\  
\text{where }\delta(x)=(x-1/x)(qx-1/(qx))
\end{eqnarray*}
This link is noted $6^3_2$ in \cite{Wit} and we check
$LG(6^3_2)=\delta(\frac p{\sqrt q})M(B)(q,\frac p{\sqrt q},\frac
p{\sqrt q},\frac p{\sqrt q})$.\\

The first link with trivial multivariable Alexander polynomial in the
Thistlethwaite link table is $L9n27$ 
$$L9n27=\epsh{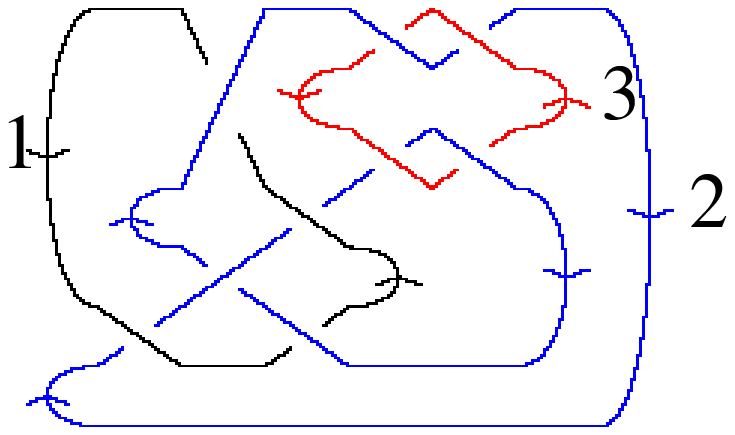}{16ex}$$ (see the ``Knot Atlas'' at the url
http://katlas.math.toronto.edu).  We have
$$M(L9n27)(q,q_1,q_2,q_3)=\qn{1}\qn{2}(q^2q_2^2+q_2^{-2}-2)$$
this is coherent with $\nabla(L9n27)=0$ and this show in particular,
that $M$ is strictly stronger than the Conway potential function.
\linespread{1}

\vfill

\end{document}